\documentclass[11pt,twoside]{article}
\usepackage{amsmath,amsfonts,amssymb,amsthm,graphicx,epsf,epstopdf,
color,pifont,subfigure,flafter,bm,amscd}
\usepackage{tikz}
\allowdisplaybreaks[4]
\usepackage{subfigure}
\usepackage{caption}

\newtheorem{lm}{Lemma}[section]
\newtheorem{thm}{Theorem}[section]

\newcounter{saveeqn}%

\makeatletter
\evensidemargin
\oddsidemargin
\makeatother

\title{\Large\bf Bifurcations and explicit unfoldings of grazing loops connecting one high multiplicity tangent point
\thanks{
Supported by NSFC \#12271378 and the Natural Science Foundation of Sichuan Province of China.
}
}

\author{Zhihao Fang$^{1}$,~~ Xingwu Chen$^{1,2}$
\footnote{Corresponding author. Email address:
scuxchen@scu.edu.cn, xingwu.chen@hotmail.com (X. Chen).}
\\
{\small 1 School of Mathematics, Sichuan University,}
{\small Chengdu, Sichuan 610065, P. R. China}\\
{\small 2 School of Mathematics and Big Data, Chongqing University of Arts and Sciences,}\\
{\small Chongqing 402160, P. R. China}
}
\date{}

\begin{document}
\maketitle

%%%%%%%%%%%%%%%%%%%%%%%%%%%%%%%%%%%%%%%%%%%%%%%%%%%%%%%%%%%%%%%%%%%%%%%

\begin{abstract}
For piecewise-smooth differential systems, in this paper we focus on crossing limit cycles and sliding loops
bifurcating from a grazing loop connecting one high multiplicity tangent point.
For the low multiplicity cases considered in previous publications, the method is to define and analyze
return maps following the classic idea of Poincar\'e. However, high multiplicity leads to
that either domains or properties of return maps are unclear under perturbations.
To overcome these difficulties, we unfold grazing loops by functional parameters and functional functions,
and analyze this unfolding along some specific parameter curve. Relationships between multiplicity
and the numbers of crossing limit cycles and sliding loops are given, and our results not only
generalize the results obtained in [J. Differential Equations 255(2013), 4403-4436; 269(2020), 11396-11434],
but also are new for some specific grazing loops.

\vskip 0.2cm


\end{abstract}

\baselineskip 15pt
\parskip 10pt
\thispagestyle{empty}
\setcounter{page}{1}

%%%%%%%%%%%%%%%%%%%%%%%%%%%%%%%%%%%%%%%%%%%%%%%%%%%%%%%%%%%%
\section{Introduction}

\setcounter{equation}{0}
\setcounter{lm}{0}
\setcounter{thm}{0}
\setcounter{rmk}{0}
\setcounter{df}{0}
\setcounter{cor}{0}
\setcounter{prop}{0}

Consider a piecewise-smooth differential system defined on a bounded open set ${\cal U}\subset\mathbb{R}^2$ containing the origin $O:(0,0)$
\begin{equation}
\begin{aligned}
		\left( \begin{array}{c}
		\dot{x}\\
		\dot{y}\\
	\end{array} \right) =\left\{
    \begin{aligned}
&		\left( \begin{array}{c}
			f^+(x,y)\\
			g^+(x,y)
		\end{array} \right)   &&\mathrm{if}~(x,y)\in\Sigma^+,\\
&		\left( \begin{array}{c}
			f^-(x,y)\\
			g^-(x,y)
		\end{array} \right)   && \mathrm{if}~(x,y)\in\Sigma^-,
	\end{aligned} \right.
    \label{pws1}
\end{aligned}
\end{equation}
where $\dot x:=dx/dt, \dot y:=dy/dt$, $f^\pm,g^\pm$ are $C^{\infty}$ functions on $\mathbb{R}^2$,
$\Sigma^+:=\left\{(x,y)\in{\cal U}:~y>0\right\}$ and $\Sigma^-:=\left\{(x,y)\in{\cal U}:~y<0\right\}$.
As defined in \cite{Bernardo08,Filippov88},
$\Sigma:=\left\{(x,y)\in{\cal U}:~y=0\right\}$ is called a {\it switching manifold}
(or {\it switching boundary})
which divides system~\eqref{pws1}
into the {\it upper subsystem} defined on $\Sigma^+$ and the {\it lower subsystem} defined on $\Sigma^-$. A continuous function $(x(t),y(t))^\top$ defined over some interval $I$ is called a {\it solution} of system~\eqref{pws1}
if it satisfies differential inclusion
$
(\dot x,\dot y)^\top\in F(x,y)
$
almost everywhere over $I$, where
\begin{equation*}
F(x,y):=\left\{
\begin{aligned}
&\left\{{\cal Z}^+(x,y)\right\},  &&(x,y)\in\Sigma^+, \\
&\left\{a{\cal Z}^+(x,y)+(1-a){\cal Z}^-(x,y):~a\in[0,1]\right\}, &&(x,y)\in\Sigma, \\
&\left\{{\cal Z}^-(x,y)\right\},  &&(x,y)\in\Sigma^-
\end{aligned}
\right.
\end{equation*}
and ${\cal Z}^\pm(x,y):=\left(f^\pm(x,y),g^\pm(x,y)\right)^\top$.

Different from classic smooth differential systems, the dynamical behaviors around $\Sigma$ are most important for system~\eqref{pws1}. An equilibrium of a subsystem is called a {\it standard equilibrium} (resp. {\it boundary equilibrium})
of system~\eqref{pws1} if it lies outside of $\Sigma$ (resp. on $\Sigma$). For $p:(x,0)\in\Sigma$ satisfying $h(x):=g^+(x,0)g^-(x,0)>0$, the orbit of one subsystem reaching $p$ connects the orbit of the other subsystem escaping from $p$, which means that there is an orbit of system~\eqref{pws1} crossing $\Sigma$ at $p$. Thus, set $\Sigma_c:=\left\{(x,0)\in\Sigma:~h(x)>0\right\}$ is called a {\it crossing region}. As indicated in \cite{Filippov88}, set $\Sigma_s:=\left\{(x,0)\in\Sigma:~h(x)<0\right\}$ is called a {\it sliding region} and a {\it sliding vector field}
\begin{equation}
X_s(x,0):= \left(\frac{f^+(x,0)g^-(x,0)-f^-(x,0)g^+(x,0)}{g^-(x,0)-g^+(x,0)},0 \right)^\top.
\label{SVF}
\end{equation}
is defined on it. Thus, the orbit of one subsystem reaching $p\in\Sigma_s$ connects a sliding orbit.

\begin{figure}[h]
\centering
\subfigure[visible]
 {
  \scalebox{0.34}[0.34]{
   \includegraphics{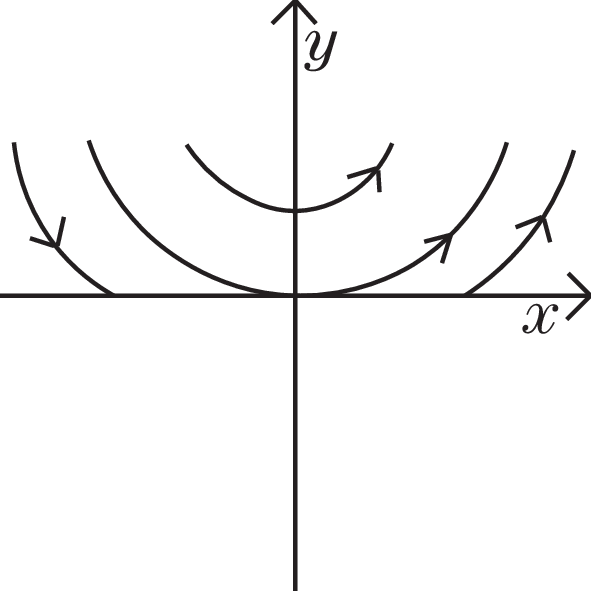}}}
\subfigure[invisible]
 {
  \scalebox{0.34}[0.34]{
   \includegraphics{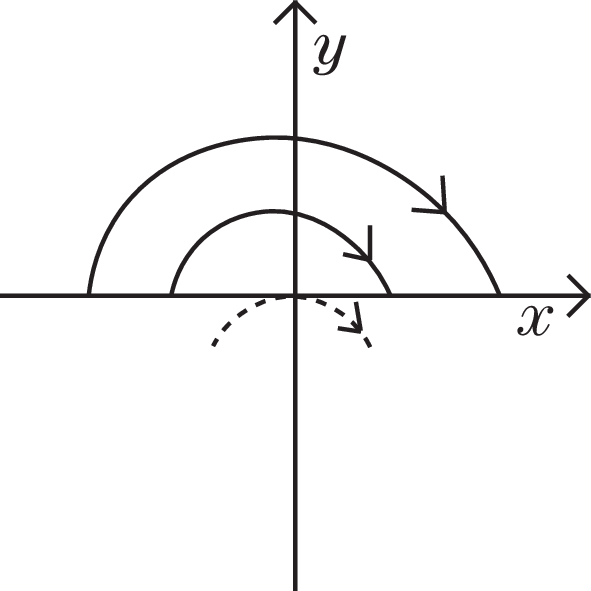}}}
\subfigure[left]
 {
  \scalebox{0.34}[0.34]{
   \includegraphics{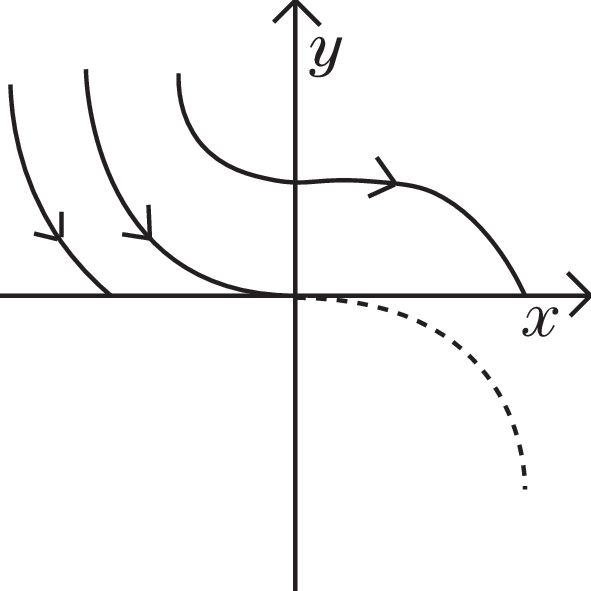}}}
\subfigure[right]
 {
  \scalebox{0.34}[0.34]{
   \includegraphics{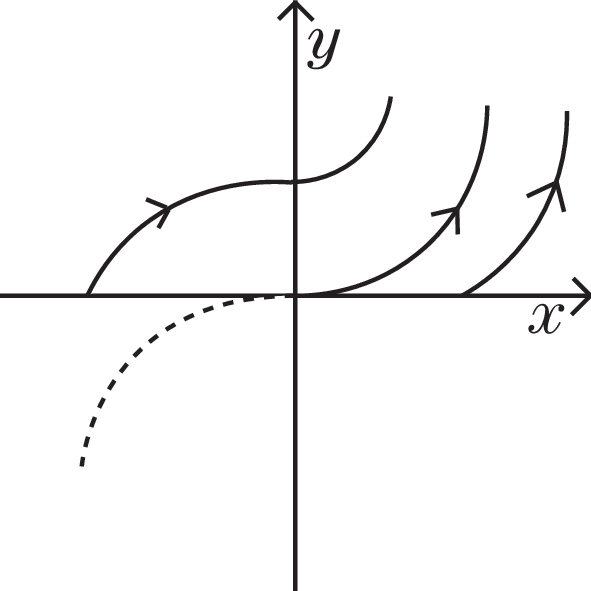}}}
   \caption{visibility of tangent point}
\label{Fig-TP}
\end{figure}

Point $(x_0,0)\in\Sigma$ is called a {\it tangent point} of system~\eqref{pws1} if $h(x_0)=0$ and $f^\pm(x_0,0)\ne0$. Clearly, tangent point implies that there is an orbit of a subsystem tangent with $\Sigma$ and does not appear in classic smooth systems. As indicated in \cite{Fang23}, tangent point $p:(x_0,0)$ is called to be of multiplicity $(m^+,m^-)$ if
\begin{eqnarray*}
&&g^+(x_0,0)=\frac{\partial g^+}{\partial x}(x_0,0)=...=\frac{\partial^{(m^+-1)}g^+}{\partial x^{(m^+-1)}}(x_0,0)=0,~\frac{\partial^{m^+}g^+}{\partial x^{m^+}}(x_0,0)\ne 0,\\
&&g^-(x_0,0)=\frac{\partial g^-}{\partial x}(x_0,0)=...=\frac{\partial^{(m^--1)}g^-}{\partial x^{(m^--1)}}(x_0,0)=0,~\frac{\partial^{m^-}g^-}{\partial x^{m^-}}(x_0,0)\ne 0,
\end{eqnarray*}
where non-negative integers $m^\pm$ satisfy $m^++m^-\ge 1$. For tangent points of the upper subsystem, we divide them into {\it visible}, {\it invisible}, {\it left} and {\it right} tangent points by location of tangent orbit as shown in Figure~\ref{Fig-TP}.
Corresponding definitions can be given for tangent points of the lower subsystem.
Further, a tangent point of system~\eqref{pws1} is called a $VI$ tangent point if it is a
visible one of the upper subsystem and an invisible one of the lower subsystem. Similarly,
other $15$ kinds $VV$, $VL$, $VR$, $IV$, $II$, $IL$, $IR$, $LV$, $LI$, $LL$, $LR$, $RV$, $RI$, $RL$, $RR$ tangent points are defined and a complete classification of tangent points is obtained.
Clearly, tangent point $p$ is of odd (resp. even) multiplicity if and only if
$p$ is either visible or invisible (resp. either left or right).

For tangent points of multiplicity $(1,0)$, it is proved
in \cite{Ponce15,Han13} that they are structurally stable,
i.e., there is no bifurcations happening under perturbations.
For tangent points of multiplicity $(0,2)$, it is proved in \cite{Kuznetsov03}
that it breaks into two tangent points of multiplicity $(0,1)$ under suitable
perturbations. According to the bifurcation phenomena, seven kinds of
tangent points of multiplicity $(1,1)$ are
given in \cite{Kuznetsov03} and another kind is given in \cite{Siller}.
For each of these eight kinds, two tangent points
of multiplicities $(0,1)$ and $(1,0)$
appear under suitable perturbations. Moreover, for one of these eight kinds it is proved
in \cite{Bonet18,Filippov88,Kuznetsov03,Han12} that there is at most
one crossing limit cycle appearing under some non-degenerate
condition and in \cite{Ponce22,Han12} that for any given integer $k$
there exist perturbations with exactly $k$ crossing limit
cycles. In \cite{Fang21,Teixeira11},
two tangent points
of multiplicity $(0,1)$ and one tangent point of multiplicity $(1,0)$ are obtained by perturbation
from a tangent point of multiplicity $(1,2)$.
For tangent points of general multiplicity $(m^+,m^-)$,
it is proved in \cite{Fang23} that for each $\ell_1\in\{1,...,m^++m^-\}$ there
exist perturbations such that there are exactly $\ell_1$ bifurcating
tangent points, i.e., tangent points bifurcating from the original tangent
point. Moreover, it is also proved in \cite{Fang23} that for an
invisible (resp. a visible, left or right) tangent point satisfying $m^+>1$
(resp. $m^+\ge 1$) of the upper subsystem and each $\ell_2\in\{1,...,(m^+-1)/2\}$
(resp. $\ell_2\in\{1,...,\lfloor(m^+ +1)/2\rfloor\}$), there exist
perturbations such that there are $\lfloor(m^+-1)/(2\ell_2)\rfloor$
(resp. $\lfloor(m^++1)/(2\ell_2)\rfloor$) orbits passing through
$\ell_2$ bifurcating tangent points. Here $\lfloor m \rfloor$ denotes the maximal integer no greater than $m$.
For the lower subsystem, there is a similar result.

It is well known that besides singular points, dynamical behaviors with
respect to return are very important for classic smooth system, e.g.,
limit cycles and homoclinic loops. It is no doubt that such dynamical behaviors
are also important for piecewise-smooth system~\eqref{pws1} and there also
exist some similar conceptions of loops. However, dynamical behaviors and
analysis are more complicated because there may exist tangent points or sliding motions on loops.

As indicated in \cite{Fang23}, Let $\Gamma_i(\alpha)$ denote a segment of
a regular orbit of a subsystem and be parameterized by $\alpha$,
where $\alpha\in[\alpha_{i1},\alpha_{i2}],~\alpha_{i1}<\alpha_{i2}$,
$i=1,2,...,n$. Further, each $\Gamma_i(\alpha)$ only intersects with
$\Sigma$ at endpoints $p_i$ and $p_{i+1}$, i.e.,
\begin{equation*}
\Gamma_i(\alpha_{i1})=p_{i}\in\Sigma,~~\Gamma_i(\alpha_{i2})=p_{i+1}\in\Sigma,~~\Gamma_i(\alpha)\notin\Sigma~{\rm for}~\alpha\in(\alpha_{i1},\alpha_{i2}).
\end{equation*}
An oriented Jordan curve consisting of $\Gamma_i(\alpha)$ ($i=1,...,n$)
is called a {\it crossing periodic orbit} of system~\eqref{pws1} if
all $p_i$ ($i=1,...,n$) are crossing points. Further, a crossing
periodic orbit is called a {\it crossing limit cycle} if it is
isolated, and denoted by $L_c$. Clearly, crossing limit cycles in
system~\eqref{pws1} are similar to limit cycles in smooth systems.

\begin{figure}[h]
\centering
\subfigure[$L^{cri}$]
 {
  \scalebox{0.33}[0.33]{
   \includegraphics{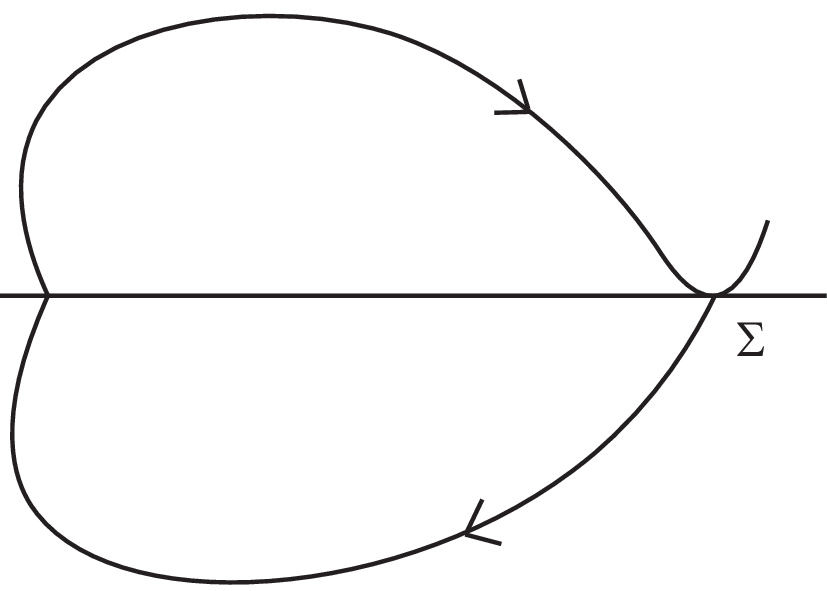}}}
\subfigure[$L^{cro}$]
 {
  \scalebox{0.33}[0.33]{
   \includegraphics{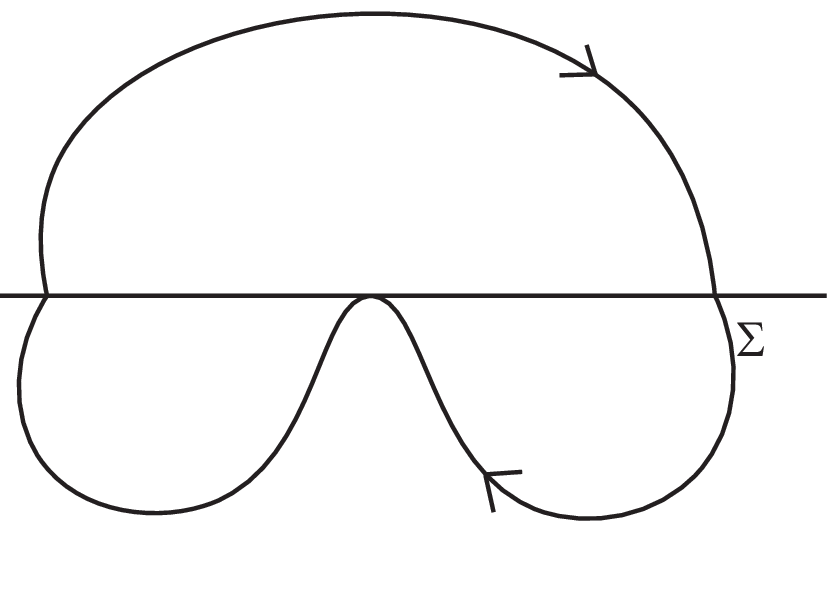}}}
\subfigure[$L^{gra}$]
 {
  \scalebox{0.33}[0.33]{
   \includegraphics{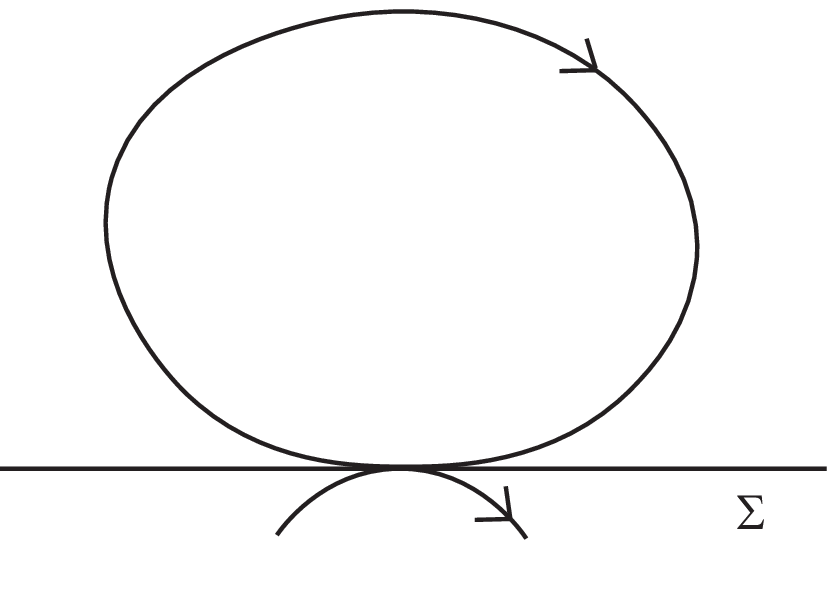}}}
   \caption{examples of $L^{cri},~L^{cro},~L^{gra}$}
\label{Fig-Loop}
\end{figure}

Different from $L_c$, an oriented Jordan curve is called a
{\it sliding loop} connecting tangent points if it consists of $\Gamma_i(\alpha)$
($i=1,...,n$) and at least one sliding orbit, and denoted by $L_s$.
As a contrary, an oriented Jordan curve is called a {\it nonsliding loop}
connecting tangent points if it consists of $\Gamma_i(\alpha)$ ($i=1,...,n$)
and there exists $i^*\in \{1,...,n\}$ such that $p_{i^*}$ is a
tangent point, and denoted by $L_{ns}$. We call point $p_i$
a {\it switching point} of $L_{ns}$ if $L_{ns}$ enters one half-plane
from the other half-plane at $p_i$. Further, $L_{ns}$ is called to be
{\it critical} (resp. {\it crossing}, {\it grazing}) if it has switching
points and some of them are tangent points (resp. it has switching points
and all of them are crossing points, it has no switching points), and
denoted by $L^{cri}$ (resp. $L^{cro}$, $L^{gra}$). Some examples of
$L^{cri}, L^{cro}, L^{gra}$ are shown in Figure~\ref{Fig-Loop}.
The bifurcations of these loops are important in dynamics and 
have extensive practical backgrounds, see e.g. \cite{Bernardo08,Bernardo082}.

For $L^{cri}$ connecting one tangent point of multiplicity $(0,1)$,
it is proved in \cite{Ponce14,Ponce15,Kuznetsov03,Han13} that there is at most
one bifurcating $L_c$ or one bifurcating $L_s$ under perturbations, and it is reachable.
For $L^{cri}$ connecting one tangent point of multiplicity $(1,1)$, it is
proved in \cite{Han13,Novaes18,Huang22} that there are either two bifurcating $L_c$
or one bifurcating $L_c$ and one bifurcating $L_s$ under perturbations.
For $L^{cri}$ connecting two tangent points of multiplicities $(1,0)$ and $(0,1)$,
it is proved in \cite{Huang22} that there is either one bifurcating $L_c$ or one
bifurcating $L_s$ under perturbations. For $L^{cri}$ connecting one tangent point
of general multiplicity $(m^+,m^-)$, it is proved in \cite{Fang23} that the
sum of the numbers of bifurcating $L_c$ and $L_s$ is at least one
and at least $\max\{m^\pm\}$ if additionally $m^\pm\ge5$.

For $L^{gra}$, by its definition the tangent points on $L^{gra}$ lie
on one periodic orbit of a subsystem, and usually of the upper subsystem without loss of generality.
Assume that this periodic orbit is a hyperbolic limit cycle. For $L^{gra}$ connecting
one tangent point of multiplicity $(1,0)$, at most one bifurcating $L_s$ exists and it is reachable as in \cite{Kuznetsov03}.
For $L^{gra}$
connecting one $VI$ tangent point of multiplicity $(1,1)$, it is proved
in \cite{Han13} that there is at least one bifurcating $L_c$. Later, it is proved in \cite{Li20} that there are
at most either two bifurcating $L_c$ or one bifurcating $L_c$ and one bifurcating
$L_s$ and they are reachable. For $L^{gra}$ connecting one $VV$
tangent point of multiplicity $(1,1)$, it is proved in \cite{Li20}
that the maximum number of bifurcating $L_s$ is one.
For $L^{cro}$, as far as we know, there is no results and one can find that the
research methods and dynamical behaviors are highly similar to $L^{gra}$.

Compared with the fact that the bifurcations of $L^{cri}$ are investigated not only for
some specific small multiplicities $m^\pm$ but also for general multiplicities,
the bifurcations of $L^{gra}$ are only focused on some specific small multiplicities $m^\pm\le1$.
To develop the bifurcations of $L^{gra}$ for general multiplicities, it is natural to firstly consider
{\it how about the bifurcations of $L^{gra}$ connecting one tangent point of multiplicity $(1,m)$ for general $m\ge2$?}

\begin{figure}[h]
\centering
\subfigure[$S$-$I$]
 {
  \scalebox{0.34}[0.34]{
   \includegraphics{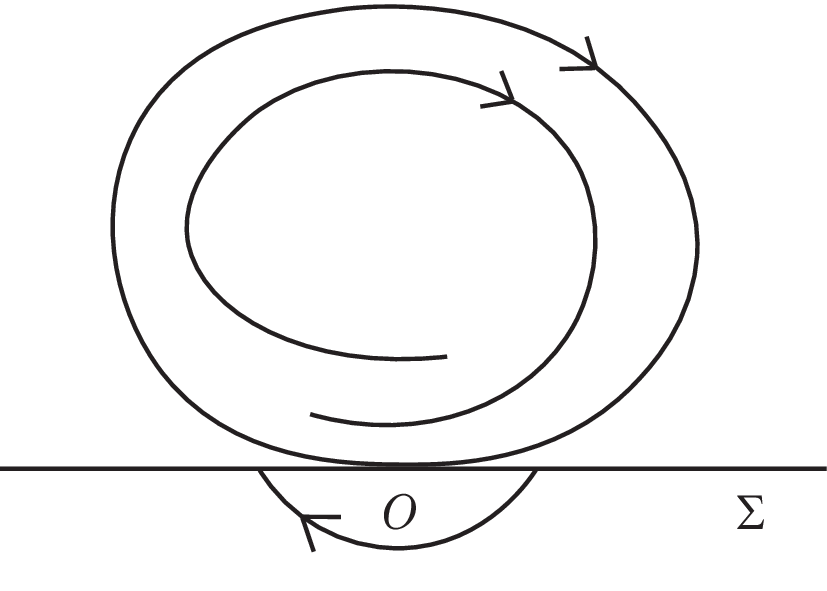}}}
\subfigure[$S$-$V$]
 {
  \scalebox{0.34}[0.34]{
   \includegraphics{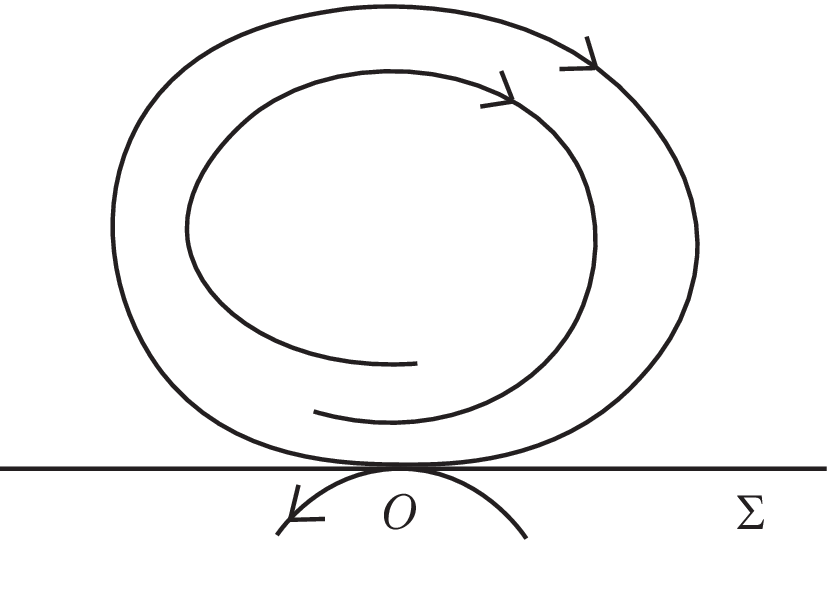}}}
\subfigure[$S$-$L$]
 {
  \scalebox{0.34}[0.34]{
   \includegraphics{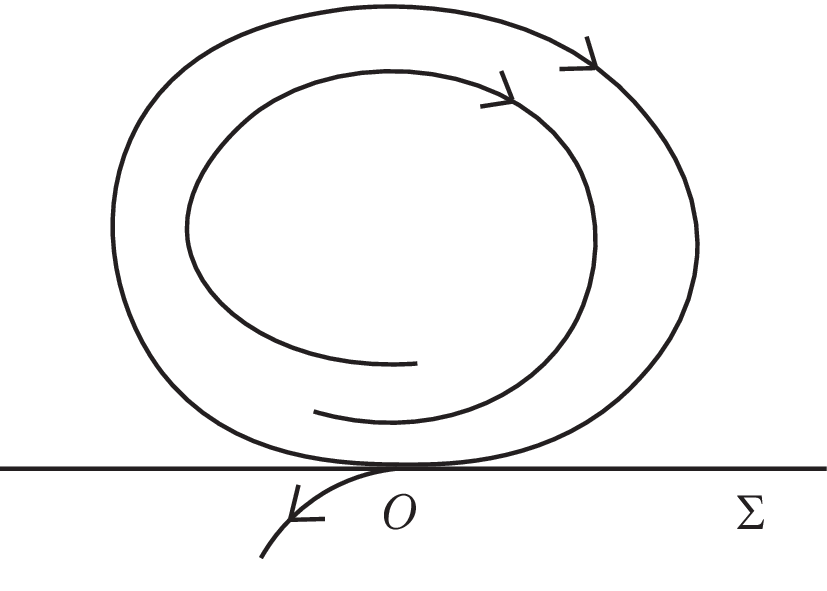}}}
\subfigure[$S$-$R$]
 {
  \scalebox{0.34}[0.34]{
   \includegraphics{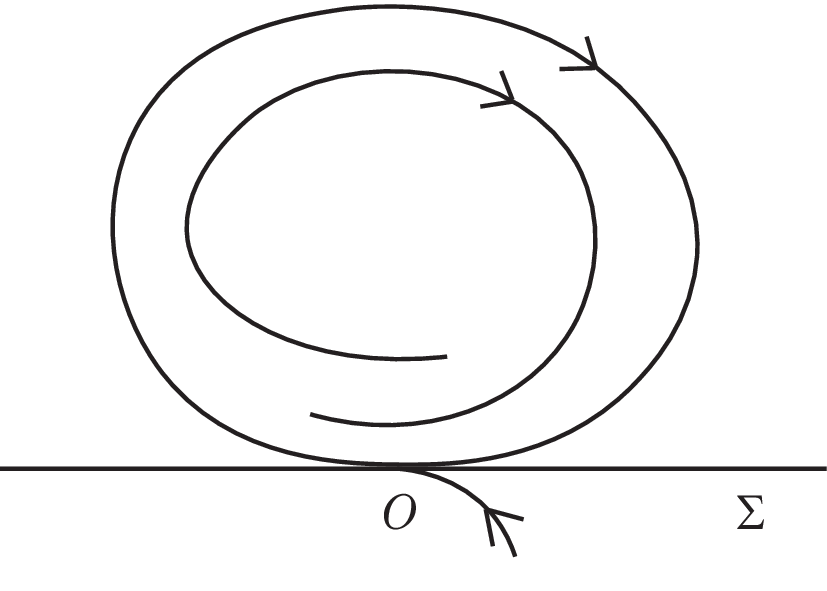}}}
\caption{ grazing loop $L^{gra}_*$ }
\label{Fig-Grazing}
\end{figure}

Motivated by this question, in this paper we investigate the numbers of crossing limit cycle $L_c$ and
sliding loop $L_s$ bifurcating from a grazing loop $L^{gra}$ connecting one tangent point of multiplicity $(1,m)$.
Let $\beta_c$ and $\beta_s$ denote the number of bifurcating $L_c$ and $L_s$ from $L^{gra}$ respectively.
It is not hard to check
in \cite{Kuznetsov03,Li20,Han13} that $\beta_c+\beta_s$ increases to two from one (resp. remains to one)
as the tangent point of multiplicity $(1,0)$
degenerates to a $VI$(resp. $VV$) tangent point of multiplicity $(1,1)$.
Thus, it is reasonable for us to claim that $\beta_c, \beta_s$ are closely related to multiplicity and
visibility of the tangent point. Then, we focus on establishing this relationship and consider system~\eqref{pws1} satisfying $g^-(x,y)=\phi(x,y)x^{m}$, i.e.,
\begin{equation}
\begin{aligned}
		\left( \begin{array}{c}
		\dot{x}\\
		\dot{y}\\
	\end{array} \right) =\left\{
    \begin{aligned}
&		\left( \begin{array}{c}
			f^+(x,y)\\
			g^+(x,y)
		\end{array} \right)   &&\mathrm{if}~(x,y)\in\Sigma^+,\\
&		\left( \begin{array}{c}
			f^-(x,y)\\
			\phi(x,y)x^m
		\end{array} \right)   && \mathrm{if}~(x,y)\in\Sigma^-,
	\end{aligned} \right.
    \label{pws2}
\end{aligned}
\end{equation}
where $m\ge 2$ and there is a grazing loop $L^{gra}_*$ connecting a unique tangent point $O$ of multiplicity $(1,m)$.
We assume that $L^{gra}_*$ is a clockwise and hyperbolic limit cycle of the upper subsystem.
Since we focus on finding bifurcating $L_c$ and $L_s$ as discussions above,
only cases satisfying $f^+(0,0)f^-(0,0)>0$ are analyzed in this paper.
There are totally eight types of $L^{gra}_*$ by the lower visibility (invisible, visible, left, right) of
tangent point $O$ and stability of this limit cycle.
Further, $L^{gra}_*$ is called to be of {\it type $S$-$I$} if this limit cycle is
stable and tangent point $O$ is invisible for the lower subsystem. Similarly,
other seven types $S$-$V$, $S$-$L$, $S$-$R$, $U$-$I$, $U$-$V$, $U$-$L$, $U$-$R$ can
be defined. Under transformation $(x,y,t)\to(-x,y,-t)$ loop $L^{gra}_*$ of
type $U$-$I$, $U$-$V$, $U$-$L$, $U$-$R$ becomes one of type $S$-$I$, $S$-$V$, $S$-$R$, $S$-$L$
respectively, which implies that we need only to consider four types as shown in Figure~\ref{Fig-Grazing}.

The main result is as follows.

\begin{thm}
Assume that $f^+(0,0)f^-(0,0)>0$ in system~\eqref{pws2}.
Then there exists a perturbation system of form
\begin{equation}
\begin{aligned}
		\left( \begin{array}{c}
		\dot{x}\\
		\dot{y}\\
	\end{array} \right) =\left\{
    \begin{aligned}
&		\left( \begin{array}{c}
			f^+(x,y+\alpha)\\
			g^+(x,y+\alpha)
		\end{array} \right)   &&\mathrm{if}~(x,y)\in\Sigma^+,\\
&		\left( \begin{array}{c}
            f^-(x,y+\psi)\\
            \phi(x,y+\psi)\prod_{i=1}^{m}(x-\lambda_i)-f^-(x,y+\psi)\dot\psi
        \end{array} \right)   && \mathrm{if}~(x,y)\in\Sigma^-
	\end{aligned} \right.
\end{aligned}
\label{pws2-unfold}
\end{equation}
such that the following statements hold.
\begin{enumerate}
       \item[{\rm (a)}] In the case that $L^{gra}_*$ is of type $S$-$I$, $S$-$L$, (resp. $S$-$V$),
       \begin{equation}
       \beta_c+\beta_s\ge m+1~~({\rm resp.~}\beta_c+\beta_s\ge m),
       \label{res-1}
       \end{equation}
        and
       \begin{equation}
       \beta_c\ge m+1-\ell~({\rm resp.~}\beta_c\ge m-\ell),~~~\beta_s=\ell
       \label{res-2}
       \end{equation}
       for each $\ell\in\left\{0,...,\lfloor m/2\rfloor \right\}$ if additionally $m\ge 4$.
       \item[{\rm (b)}] In the case that $L^{gra}_*$ is of type $S$-R,
      $
       \beta_c+\beta_s\ge 2
     $
       if $m=2$, and
       \begin{equation}
       \beta_c\ge m+2-\ell,~~~\beta_s=\ell
       \label{res-4}
       \end{equation}
       for each $\ell\in\left\{0,...,m/2\right\}$ if $m\ge4$.
\end{enumerate}
Here $\psi:=\psi(x,{\boldsymbol k})$ is a $C^{\infty}$ function satisfying $\psi(x,{\boldsymbol 0})\equiv 0$, $\dot\psi:=d\psi/dx$, $\alpha\in\mathbb{R},{\boldsymbol \lambda}:=(\lambda_1,...,\lambda_m)\in\mathbb{R}^m,~{\boldsymbol k}:=(k_1,...,k_n)\in\mathbb{R}^n$.
\label{thm1}
\end{thm}

Theorem~\ref{thm1} gives a generalization on $\beta_c, \beta_s$ from $m=1$ to general $m$ for $L^{gra}_*$ of type $S$-$I$, $S$-$V$ and gives new results for $L^{gra}_*$ of type $S$-$L$, $S$-$R$. In the later proof, suitable $\psi(x,{\boldsymbol k})$ can be defined for $m\ge4$ but invalid for $m=2,3$, which explains why \eqref{res-2} and \eqref{res-4} only hold for $m\ge4$.

This paper is organized as follows. Some necessary preliminary theories for proof are stated in section 2
and a proof of Theorem~\ref{thm1} is given in section 3. We summarize conclusions and give some discussions in section 4 to end this paper.

%%%%%%%%%%%%%%%%%%%%%%%%%%%
\section{Preliminary lemmas}

\setcounter{equation}{0}
\setcounter{lm}{0}
\setcounter{thm}{0}
\setcounter{rmk}{0}
\setcounter{df}{0}
\setcounter{cor}{0}

To investigate $\beta_c$ and $\beta_s$, the first thing is to find appropriate function  $\psi(x,{\boldsymbol k})$
such that system~\eqref{pws2-unfold} is a small perturbation of system~\eqref{pws2}, that is,
the vector field $X(x,y)$ of \eqref{pws2-unfold} lies in a small neighborhood of the vector field $Y(x,y)$ of \eqref{pws2}.
Denote corresponding upper vector fields and lower ones by
\begin{equation*}
\begin{aligned}
&X^\pm(x,y):=\left(X^\pm_1(x,y),X^\pm_2(x,y)\right)^\top\in C^\infty\left(\mathbb{R}^2,\mathbb{R}^2\right),\\
&Y^\pm(x,y):=\left(Y^\pm_1(x,y),Y^\pm_2(x,y)\right)^\top\in C^\infty\left(\mathbb{R}^2,\mathbb{R}^2\right).
\end{aligned}
\end{equation*}
As indicated in \cite{Fang23}, the distance of system~\eqref{pws2} and system~\eqref{pws2-unfold} is
\begin{eqnarray*}
\rho\left(X,Y\right):=d\left(X^+_1,Y^+_1\right)+d\left(X^+_2,Y^+_2\right)+d\left(X^-_1,Y^-_1\right)+d\left(X^-_2,Y^-_2\right),
\end{eqnarray*}
where
\begin{equation*}
d\left(P, Q\right):=\max_{(x,y)\in\bar{\cal U}}\left\{\sum_{i_1+i_2=0}^{1}\left|\frac{\partial^{i_1+i_2} (P-Q)}{\partial x^{i_1}\partial y^{i_2}}\right|\right\}
\end{equation*}
for $P(x,y),Q(x,y)\in C^\infty\left(\mathbb{R}^2,\mathbb{R}\right)$.

\begin{figure}[h]
\centering
\includegraphics[scale=0.4]{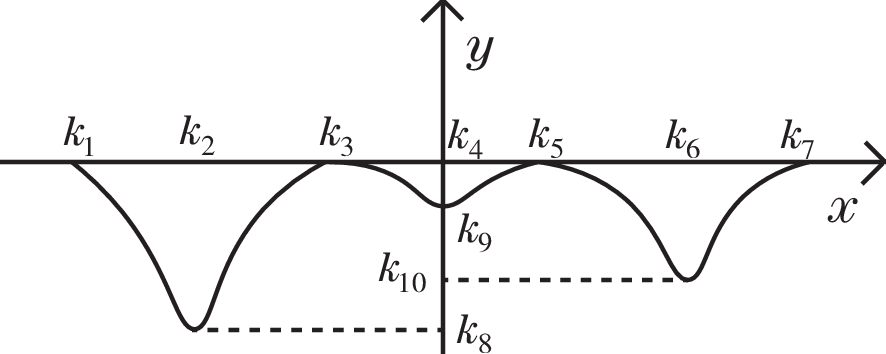}
\caption{$\psi(x,{\boldsymbol k})$ with $d=3$ and ${\boldsymbol k}\in{\cal K}, k_8,k_9,k_{10}<0$}
\label{Fig-psi}
\end{figure}

 Consider $C^\infty$ cutoff functions (see, e.g., \cite{Lee18})
\begin{equation*}
h^*(x,r_1,r_2):=\left\{
\begin{aligned}
&0,&&x\in (-\infty,r_1], \\
&\frac{1}{1+e^{\eta(x)}},&&x\in(r_1,r_2), \\
&1,&&x\in[r_2,+\infty),
\end{aligned}
\right.~~
h_*(x,r_1,r_2):=\left\{
\begin{aligned}
&1,&&x\in(-\infty,r_1], \\
&\frac{e^{\eta(x)}}{1+e^{\eta(x)}},&&x\in(r_1,r_2), \\
&0,&&x\in[r_2,+\infty),
\end{aligned}
\right.
\end{equation*}
where $r_1<r_2$ and
\begin{equation*}
\eta(x):=\frac{1}{x-r_1}+\frac{1}{x-r_2}.
\end{equation*}
Then for $x\in\mathbb{R}$ and any integer $d>0$, we define
\begin{eqnarray*}
\psi(x,{\boldsymbol k}):=
\left\{
\begin{aligned}
&\psi^*(x,{\boldsymbol k}), &&{\rm when}~{\boldsymbol k}\in{\cal K},\\
&0,~~~~~~&&{\rm when}~{\boldsymbol k}\in\mathbb{R}^{3d+1}\setminus{\cal K},
\end{aligned}
\right.
\end{eqnarray*}
where
\begin{equation*}
\psi^*(x,{\boldsymbol k}):=
\left\{
\begin{aligned}
& k_{2d+1+i}h^*\left(x,k_{2i-1},k_{2i}\right),~~~~&&x\in\left(k_{2i-1},k_{2i}\right],\\
& k_{2d+1+i}h_*\left(x,k_{2i},k_{2i+1}\right),~~~~&&x\in\left(k_{2i},k_{2i+1}\right],\\
& 0,&&{\rm otherwise}
\end{aligned}
\right.
\end{equation*}
and $i=1,...,d, {\cal K}:=\left\{{\boldsymbol k}\in{\mathbb R}^{3d+1}:~k_1< ...<k_{2d+1}\right\}$.
It is not hard to obtain that $\psi(x,{\boldsymbol k})$ is a $C^{\infty}$ function and satisfies
$\psi(x,{\boldsymbol 0})\equiv 0$. Moreover, derivatives of $\psi(x,{\boldsymbol k})$ of any order with respect to $x$ are all zero at $k_i$ ($i=1,...,2d+1$) (see, e.g., \cite{Lee18}). An example is shown in Figure~\ref{Fig-psi}.

The following lemma is proved in \cite{Fang23} and gives a sufficient condition for the vector field $X(x,y)$ of \eqref{pws2-unfold} lies in a small neighborhood of the vector field $Y(x,y)$ of \eqref{pws2}.

\begin{lm}
If ${\boldsymbol k}$ in function $\psi(x,{\boldsymbol k})$ satisfies either ${\boldsymbol k}\in \mathbb{R}^{3d+1}\setminus {\cal K}$ or ${\boldsymbol k}\in {\cal K}$ with
\begin{equation}
k_{i+2d+1}=o\left(\left(k_{2i}-k_{2i-1}\right)^4\right),~~~k_{2i+1}-k_{2i}=\beta_i(k_{2i}-k_{2i-1})
\label{K}
\end{equation}
for some constants $\beta_i>0,~i=1,...,d$,
then for any $\epsilon>0$ there exists $\delta>0$ such that
\begin{equation*}
\rho\left(X,Y\right)<\epsilon,~~~~\forall |\alpha|,|{\boldsymbol \lambda}|, |{\boldsymbol k}|<\delta.
\end{equation*}
\label{lm1}
\end{lm}

In order to analysis the perturbation system~\eqref{pws2-unfold}, we needs a so-called {\it transition system},
which is \eqref{pws2-unfold} satisfying $\alpha=0$ and ${\boldsymbol k}={\boldsymbol 0}$, i.e.,
\begin{equation}
\begin{aligned}
		\left( \begin{array}{c}
		\dot{x}\\
		\dot{y}\\
	\end{array} \right) =\left\{
    \begin{aligned}
&		\left( \begin{array}{c}
			f^+(x,y)\\
			g^+(x,y)
		\end{array} \right)   &&\mathrm{if}~(x,y)\in\Sigma^+,\\
&		\left( \begin{array}{c}
			f^-(x,y)\\
			\phi(x,y)\prod_{i=1}^{m}(x-\lambda_i)
		\end{array} \right)   && \mathrm{if}~(x,y)\in\Sigma^-.
	\end{aligned} \right.
    \label{pws2-tran}
\end{aligned}
\end{equation}
The following lemma is proved in \cite{Fang23} and gives the relation between solutions of the lower subsystems
of \eqref{pws2-unfold} and \eqref{pws2-tran}.

\begin{lm}
Assume that ${\boldsymbol k}\in\mathbb{R}^{3d+1}\setminus{\cal K}$ or ${\boldsymbol k}\in{\cal K}$ satisfying~\eqref{K}.
Then for $(x_0,y_0)\in\Sigma^-$ there exist $\delta>0$ such that for $|{\boldsymbol k}|<\delta$,
\begin{equation}
\left(
\begin{aligned}
&\widetilde\gamma^-_1\left(t,x_0,y_0-\psi\left(x_0,{\boldsymbol k}\right)\right)\\
&\widetilde\gamma^-_2\left(t,x_0,y_0-\psi\left(x_0,{\boldsymbol k}\right)\right)
\end{aligned}
\right)=
\left(
\begin{aligned}
&\widehat\gamma^-_1\left(t,x_0,y_0\right)\\
&\widehat\gamma^-_2\left(t,x_0,y_0\right)
\end{aligned}
\right)
-
\left(
\begin{aligned}
&0\\
&\psi\left(\widehat \gamma^-_1\left(t,x_0,y_0\right),{\boldsymbol k}\right)
\end{aligned}
\right)
\label{defor}
\end{equation}
over some intervals $I$, where
\begin{equation*}
\begin{aligned}
&\widetilde\gamma^-\left(t,x_0,y_0\right):=\left(\widetilde\gamma^-_1\left(t,x_0,y_0\right), \widetilde\gamma^-_2\left(t,x_0,y_0\right)\right)^\top,\\
&\widehat\gamma^-\left(t,x_0,y_0\right):=\left(\widehat\gamma^-_1\left(t,x_0,y_0\right), \widehat\gamma^-_2\left(t,x_0,y_0\right)\right)^\top
\end{aligned}
\end{equation*}
denote solutions the lower subsystem of~\eqref{pws2-unfold} and system~\eqref{pws2-tran} with initial value $(x_0,y_0)$.
\label{lm2}
\end{lm}

\begin{figure}[htp]
\centering
\includegraphics[scale=0.38]{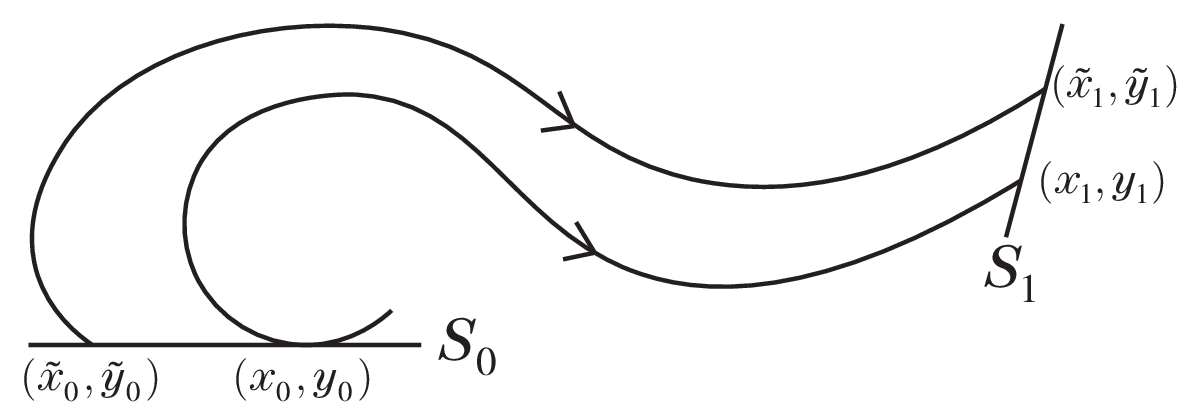}
\caption{transition map $V(r)$}
\label{Fig-trans}
\end{figure}

In the rest of this section, we introduce {\it transition map} and give some corresponding properties.
Consider the following system
\begin{equation}
        \left( \begin{array}{c}
        \dot{x}\\
        \dot{y}\\
    \end{array} \right)
    =
      \left( \begin{array}{c}
            f(x,y)\\
            g(x,y)
        \end{array} \right)
\label{trans-sys}
\end{equation}
where $(x,y)\in{\mathbb R}^2$ and $f, g$ are $C^{\infty}$. Let $\gamma(t,x_0,y_0):=\left(\gamma_1(t,x_0,y_0),\gamma_2(t,x_0,y_0)\right)^\top$
be the solution of system~\eqref{trans-sys} with initial value $(x_0,y_0)$. For a regular point $(x_0,y_0)$, we assume that orbit $\gamma\left(t,x_0,y_0\right)$ is tangent with a horizontal sector $S_0$ at $(x_0,y_0)$
but transversally intersects $S_1$ at $(x_1,y_1)$ at $t=T$ as shown in Figure~\ref{Fig-trans}. Further,
for any point $\left(\widetilde x_0,\widetilde y_0\right)$ sufficiently close to $(x_0,y_0)$ on
$S_0$, orbit $\gamma\left(t,\widetilde x_0,\widetilde y_0\right)$ intersects $S_1$ at $(\widetilde x_1,\widetilde y_1)$. Let $N_0:=(N_{01},0)^\top$ (resp. $N_1:=(N_{11},N_{12})^\top$) be the unit vector parallel to $S_0$ (resp. $S_1$), we write $(\widetilde x_0,\widetilde y_0)=\left(x_0+rN_{01}, y_0\right)$ and
\begin{equation*}
\left(\widetilde x_1,\widetilde y_1\right)=\left(x_1+V(r)N_{11}, y_1+V(r)N_{12}\right)
\end{equation*}
for sufficiently small $r$. As indicated in \cite{JKHaleBook,KuznetsovBook,ZhangJ98} function $V(r)$ is called a {\it transition map} and it is characterized by the following lemma, which is proved in \cite{Fang23}.

\begin{lm}
$V(r)$ is $C^{\infty}$ and if
\begin{equation*}
g(x_0,y_0)=0,~~~~~~\frac{\partial g}{\partial x}(x_0,y_0)\ne0,
\end{equation*}
then
\begin{equation*}
V\left(r\right)=V_{2}r^{2}+O\left(r^{3}\right),
\end{equation*}
where
\begin{equation*}
V_{2}:=\frac{\partial g}{\partial x}(x_0,y_0)\frac{1}{2\Delta_1}\exp\left\{\int_0^{T}\frac{\partial f\left(\gamma\left(s,x_0,y_0\right)\right)}{\partial x}
+\frac{\partial g\left(\gamma\left(s,x_0,y_0\right)\right)}{\partial y}ds\right\}\ne0
\end{equation*}
and
\begin{equation*}
\Delta_1:=\det\left(\left(f(x_0,y_0),g(x_0,y_0)\right)^\top,N_1\right).
\end{equation*}
\label{lm4}
\end{lm}

%%%%%%%%%%%%%%%%%%%%%%%%%%%%%%%%%%%%%%%%%%%%%%%%%%%%%%%%%%%
\section{Proof of Theorem~\ref{thm1}}
\setcounter{equation}{0}
\setcounter{lm}{0}
\setcounter{thm}{0}
\setcounter{rmk}{0}
\setcounter{df}{0}
\setcounter{cor}{0}

In this section, with those lemmas given in last section we give a proof for our main result.

\begin{proof}[Proof of Theorem~\ref{thm1}(a)]
We begin with the case that $L^{gra}_*$ of type $S$-$I$.
For the lower subsystem of~\eqref{pws2-tran}, it is
not difficult to obtain that $m$ is odd. Thus, we need only to consider two cases: odd $m\ge 5$ and $m=3$,
and split the proof into three steps.

{\it Step 1. Analyze the perturbation of the upper subsystem via transition maps}.
Consider the upper subsystem of transition system~\eqref{pws2-tran}. Take a point $(a,b)$ in $L^{gra}_*$ satisfying $g^+(a,b)\ne0,~a>0$ and three sectors
\begin{equation*}
S_1:=\left\{(x,0):~x\in(-\epsilon,0]\right\},~S_2:=\left\{(x,0):~x\in[0,\epsilon)\right\},~S_3:=\left\{(x,b):~x\in(a-\epsilon,a+\epsilon)\right\}
\end{equation*}
for $0<\epsilon\ll1$ as shown in Figure~\ref{Fig-SI-1}.
\begin{figure}[h]
\centering
\includegraphics[scale=0.5]{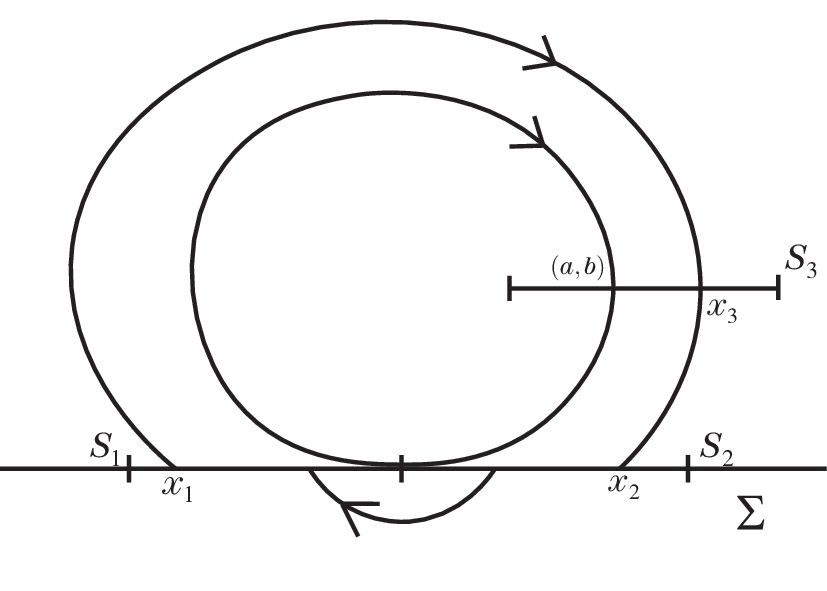}
\caption{sectors $S_1, S_2, S_3$}
\label{Fig-SI-1}
\end{figure}Let $\widehat\gamma^+(t,x_0,y_0):=\left(\widehat\gamma^+_1(t,x_0,y_0),\widehat\gamma^+_2(t,x_0,y_0)\right)^\top$
denote the solution of the upper subsystem of~\eqref{pws2-tran} with initial value $(x_0,y_0)$. Thus, there exists $T_1>0$ and $T_2<0$ such that $\widehat\gamma^+(T_1,0,0)=\widehat\gamma^+(T_2,0,0)=(a,b)^\top.$
For any point $(x_1,0)\in S_1$, we obtain that orbit $\widehat\gamma^+\left(t,x_1,0\right)$ intersects $S_2$ and $S_3$ at $(x_2,0)$ and $(x_3,b)$ respectively, which means that we can define transition maps
$V^+_1(x_1)$, $V^+_2(x_2)$ via the orbits from $S_1$ to $S_3$, from $S_2$ to $S_3$, respectively.
By Lemma~\ref{lm4}, there exists $V^+_{12}>0$ and $V^+_{22}>0$ such that
\begin{equation*}
V^+_1(x_1)=a+V^+_{12}x^2_1+O(x^3_1),~~~~V^+_2(x_2)=a+V^+_{22}x^2_2+O(x^3_2).
\end{equation*}
Thus, when $V^+_1(x_1)=V^+_2(x_2)$, we get
\begin{equation}
\frac{|x_1|}{|x_2|}=\sqrt{\frac{V^+_{22}+O(x_2)}{V^+_{12}+O(x_1)}}
\label{Rela}
\end{equation}
for any sufficiently small $x_2\ne0$ and
\begin{equation*}
\sqrt\frac{V^+_{22}}{V^+_{12}}=\exp\left\{-\frac{1}{2}\int_{T_2}^{T_1}\frac{\partial f\left(\gamma\left(s,x_0,y_0\right)\right)}{\partial x}+\frac{\partial g\left(\gamma\left(s,x_0,y_0\right)\right)}{\partial y}ds\right\}>1.
\end{equation*}

{\it Step 2. Consider the case that $m\ge 5$}.
In such case, we take
\begin{equation}
\lambda_i=(i-1)\delta,~~i=1,...,m
\label{def-lam}
\end{equation}
for $0<\delta\ll1$ such that all $(\lambda_i,0)\in S_2$ and then orbits $\widehat\gamma^+\left(t,\lambda_i,0\right), i=2,...,m$ intersect sector $S_1$ at points $(\lambda^-_i,0)$ respectively for backward direction. For simplification, we rewrite points $(\lambda^-_m,0),...,(\lambda^-_2,0),(\lambda_1,0),...,(\lambda_m,0)$ as $(A_n,0), n=1,...,2m-1$ respectively. Further, it is not difficult to obtain that points $(A_n,0)\in\Sigma_c$ for $n=1,...,m-1$, points $(A_n,0)$ are invisible tangent points for $n=m, m+2,...,2m-1$ and points $(A_n,0)$ are visible tangent points for $n=m+1, m+3,...,2m-2$. Moreover, by \eqref{Rela} there exist some $M_n>0, n=1,2,...,m-1$ such that
\begin{equation}
|A_n|=\sqrt\frac{V^+_{22}+O(A_{2m-n})}{V^+_{12}+O(A_n)}A_{2m-n}=M_n\delta.
\label{Rela-del}
\end{equation}
Then for $n=2, 4,..., m-1$, we take sectors $S(A_n):=\left\{(A_n,y):~y\in(-1,1)\right\}$ and consider functions $y=\Theta_n(x)$ satisfying the Cauchy problems respectively
\begin{equation*}
\frac{dy}{dx}=\frac{\phi(x,y)\prod_{i=1}^{m}(x-\lambda_i)}{f^-(x,y)},~~y(A_1)=-\delta^m(1+n\delta):=C_n.
\end{equation*}
For all $x\in[A_1,A_{2m-1}]$, we obtain
\begin{equation}
\begin{aligned}
\Theta_n(x)
& = \Theta_n\left(A_1\right)+\int_{A_1}^{x}\frac{\phi(s,\Theta_n(s))\prod_{i=1}^{m}(s-\lambda_i)}{f^-(s,\Theta_n(s))}ds\\
& = -\delta^{m}(1+n\delta)+\left(x-A_1\right)\frac{\phi(\sigma_n,\Theta_n(\sigma_n))\prod_{i=1}^{m}(\sigma_n-\lambda_i)}{f^-(\sigma_n,\Theta_n(\sigma_n))}\\
& = -\delta^{m}(1+n\delta)+\delta^{m+1}(x^*+M_1)\frac{\phi\left(\sigma^*_n\delta,\Theta_n(\sigma^*_n\delta)\right)}{f^-(\sigma^*_n\delta,\Theta_n(\sigma^*_n\delta))}\prod_{i=1}^{m}(\sigma^*_n-i-1)\\
& = -\delta^{m}+O\left(\delta^{m+1}\right)\\
& < 0
\end{aligned}
\label{Theta-Esti}
\end{equation}
for sufficiently small $\delta>0$, where $x^*=x/\delta\in[-M_1,m-1]$, $\sigma_n\in\left(A_1,x\right)$ and $\sigma^*_n=\sigma_n/\delta\in[-M_1,x^*]$. Further, for each $n\in\left\{2, 4,...,m-1\right\}$ orbit $\widehat\gamma^-(t,A_1,C_n)$ intersects with sectors $S(A_l)=\left\{(A_l,y):~y\in(-1,1)\right\}, l=2, 4,..., m-1$ at points $(A_l,\Theta_n(A_l))$ respectively.

Then we take $d=m-1$ and
\begin{equation*}
\left\{
\begin{aligned}
&k_i=A_i,&&{\rm~for~}i=1,...,2m-1,\\
&k_{i+2d+1}=\Theta_{m+1-2i}(A_{2i}),&&{\rm~for~}i=1,...,(m-1)/2,\\
&k_{i+2d+1}=\Theta_{2i-m+1}(A_{2i}),&&{\rm~for~}i=(m+1)/2,...,m-1
\end{aligned}
\right.
\end{equation*}
in $\psi(x,{\boldsymbol k})$.
By \eqref{defor} in Lemma~\ref{lm2}, system~\eqref{pws2-unfold} with $\alpha=0$ has a grazing loop connecting origin $O$ (denoted by $L^{gra}(A_m)$) and $(m-1)/2$ critical loops connecting $(A_n,0), n=m+1, m+3...,2m-2$ (denoted by $L^{cri}(A_n)$) respectively as shown in Figure~\ref{Fig-SI-2} for $m=5$. Location of all loops satisfies
\begin{equation}
L^{gra}(A_m)~\hookrightarrow~L^{cri}(A_{m+1})~\hookrightarrow~...~\hookrightarrow~L^{cri}(A_{2m-2}),
\label{Loca}
\end{equation}
where $L^{gra}(A_m)~\hookrightarrow~L^{cri}(A_{m+1})$ means that $L^{gra}(A_m)$ lies in the region surrounded by $L^{cri}(A_{m+1})$.
\begin{figure}[h]
\centering
\includegraphics[scale=0.38]{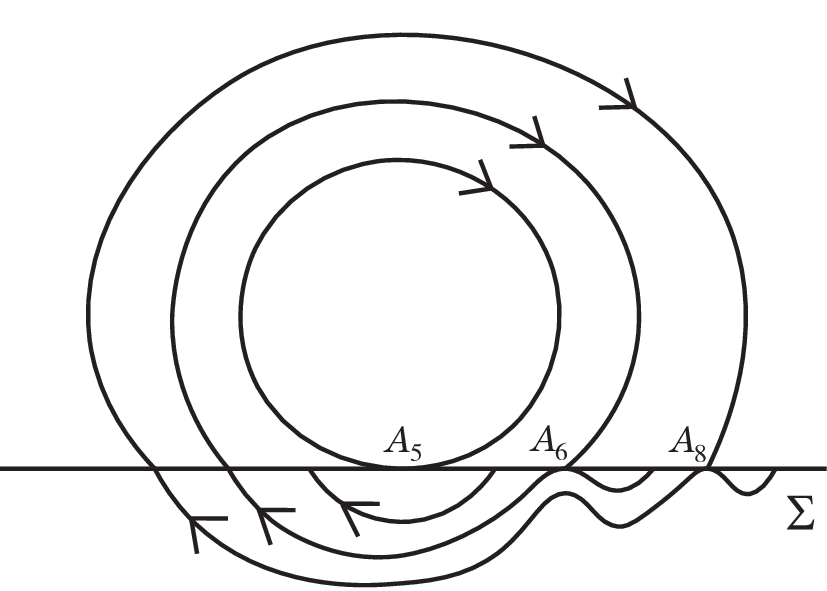}
\caption{critical loops for $m=5$}
\label{Fig-SI-2}
\end{figure}

Then we perturb grazing loop $L^{gra}(A_m)$ and those critical loops $L^{cri}(A_n)$ to obtain bifurcating $L_c, L_s$ for three cases
\begin{equation*}
{\rm~(C1)~}\ell=0,~~~~~~~{\rm~(C2)~}\ell=(m-1)/2,~~~~~~~{\rm~(C3)~}1\le\ell\le(m-3)/2.
\end{equation*}
For case (C1), as indicated in \cite{Fang23}, for any $n\in\left\{m+1, m+3,...,2m-2\right\}$ loop $L^{cri}(A_n)$ is stable, i.e., there exists point $({\cal P}_n,0)$ satisfying $0<A_n-{\cal P}_n\ll1$ such that orbit starting from it returns to $\Sigma$ at $({\cal P}^*_n,0)$ which satisfies
\begin{equation}
{\cal P}^*_n>{\cal P}_n.
\label{inner-boun}
\end{equation}
We keep $k_i, i=1,...,2m-1$ and $k_{i+2d+1}, i=1,...,(m-1)/2$ fixed and take
\begin{equation*}
k_{i+2d+1}=\Theta_{2i-m+1}(A_{2i})+a_j,~j:=i-(m-1)/2,~i=(m+1)/2,...,m-1
\end{equation*}
for sufficiently small $a_j>0$. Further, inequalities~\eqref{inner-boun} still hold and by Lemma~\ref{lm2}, for any $n\in\left\{m+1, m+3,...,2m-2\right\}$ orbit starting from $(A_n,0)$ returns to $\Sigma$ at $(A^*_n,0)$ which satisfies
\begin{equation}
A^*_n<A_n.
\label{outter-boun}
\end{equation}
Thus there is at least one bifurcating $L_c$ in a small neighborhood of $L^{cri}(A_n)$, i.e., at least $(m-1)/2$ bifurcating $L_c$ appear. Then it is proved by \cite{Li20} that there exists sufficiently small $\alpha\ne0$ and $\lambda_1\ne0$ such that there is one unstable bifurcating $L_c$ (denoted by $L^u_c$) and one stable bifurcating $L_c$ (denoted by $L^s_c$) in a small neighborhood of grazing loop $L^{gra}(A_m)$. Meanwhile, inequalities~\eqref{inner-boun} and \eqref{outter-boun} still hold, which implies that there are at least $2+(m-1)/2=(m+3)/2$ bifurcating $L_c$. Let $L^n_c$ denote the $L_c$ bifurcating from loop $L^{cri}(A_n)$ and by \eqref{Loca}, we obtain
\begin{equation}
L^u_c~\hookrightarrow~L^s_c~\hookrightarrow~L^{m+1}_c~\hookrightarrow~...~\hookrightarrow~L^{2m-2}_c.
\label{Loca-1}
\end{equation}
Since $L^s_c$ is stable and \eqref{inner-boun} holds for $n=m+1$, there is at least one bifurcating $L_c$ in the region surrounded by $L^s_c$ and $L^{m+1}_c$. On the other hand, orbit starting from $(A_n,0)$ intersects $\Sigma$ at $(A^-_n,0)$ for backward direction and then there exists $({\cal Q}_n,0)$ satisfying $0<{\cal Q}_n-A^-_n\ll1$ such that orbit starting from it returns to $\Sigma$ at $({\cal Q}^*_n,0)$ which satisfies
\begin{equation}
{\cal Q}^*_n<{\cal Q}_n.
\label{inner-boun2}
\end{equation}
Further by \eqref{inner-boun} and \eqref{inner-boun2}, there is at least one bifurcating $L_c$ in the region surrounded by $L^n_c$ and $L^{n+2}_c$ for any $n\in\left\{m+1,...,2m-4\right\}$, i.e., there are another $(m-1)/2$ bifurcating $L_c$ besides the $(m+3)/2$ bifurcating $L_c$ in \eqref{Loca-1}. Therefore, there are at least $(m+3)/2+(m-1)/2=m+1$ bifurcating $L_c$ and no bifurcating $L_s$.

For case (C2), we keep $k_i, i=1,...,2m-1$ and $k_{i+2d+1}, i=1,...,(m-1)/2$ fixed and take
\begin{equation*}
k_{i+2d+1}=\Theta_{2i-m+1}(A_{2i})-a_j,~j:=i-(m-1)/2,~i=(m+1)/2,...,m-1
\end{equation*}
for sufficiently small $a_j>0$. Further, it is proved by Lemma~\ref{lm2} that for any $n\in\left\{m+1, m+3,...,2m-2\right\}$ loop $L^{cri}(A_n)$ breaks and orbit starting from $(A_n,0)$ returns to $\Sigma$ at $(A^*_n,0)$ which satisfies
\begin{equation}
A^*_n>A_n.
\label{outter-boun2}
\end{equation}
Further, it is not difficult to obtain $(A^*_n,0)\in\Sigma_s$ and the first component of sliding vector field $X_s(x,0)$ in \eqref{SVF} satisfies
\begin{equation*}
\frac{f^+(A^*_n,0)g^-(A^*_n,0)-f^-(A^*_n,0)g^+(A^*_n,0)}{g^-(A^*_n,0)-g^+(A^*_n,0)}<0,
\end{equation*}
which implies that such orbit returns to $(A_n,0)$ along the sliding region, i.e., there is one bifurcating $L_s$ in a small neighborhood of $L^{cri}(A_n)$ (denoted by $L^n_s$). Then it is proved by a similar way that there exists sufficiently small $\alpha\ne0$ and $\lambda_1\ne0$ such that corresponding $L^u_c, L^s_c$ appear, all the $(m-1)/2$ bifurcating $L^n_s$ remain and
\begin{equation}
L^u_c~\hookrightarrow~L^s_c~\hookrightarrow~L^{m+1}_s~\hookrightarrow~...~\hookrightarrow~L^{2m-2}_s.
\label{Loca-2}
\end{equation}
Further, there is at least one bifurcating $L_c$ in the region surrounded by $L^s_c$ and $L^{m+1}_s$ because $L^s_c$ is stable and \eqref{outter-boun2} holds for $n=m+1$. Then by inequalities \eqref{inner-boun2} and \eqref{outter-boun2}, there is at least one bifurcating $L_c$ in the region surrounded by $L^n_s$ and $L^{n+2}_s$ for any $n\in\left\{m+1,...,2m-4\right\}$, i.e., there are another $(m-1)/2$ bifurcating $L_c$ besides the $2$ bifurcating $L_c$ in \eqref{Loca-2}. Therefore, there are at
least $2+(m-1)/2=(m+3)/2$ bifurcating $L_c$ and exactly $(m-1)/2$ bifurcating $L_s$.

For case (C3), it can be proved by a similar way that there exists ${\boldsymbol k}$ and $\alpha$ such that
there are at least $(m+3)/2-\ell$ bifurcating $L_c$ and exactly $\ell$ bifurcating $L_s$, which satisfy
\begin{equation}
L^u_c~\hookrightarrow~L^s_c~\hookrightarrow~L^{m+1}_c~\hookrightarrow~...~\hookrightarrow~L^{2m-2\ell-2}_c~\hookrightarrow~L^{2m-2\ell}_s~\hookrightarrow~...~\hookrightarrow~L^{2m-2}_s.
\label{Loca-3}
\end{equation}
Further, by analysis in case (C1) there is at least $(m-2\ell-1)/2$ bifurcating $L_c$ in the region surrounded by $L^s_c$ and $L^{2m-2\ell-2}_c$. On the other hand, there is at least $\ell-1$ bifurcating $L_c$ in the region surrounded by $L^{2m-2\ell}_s$ and $L^{2m-2}_s$ via analysis in case (C2). Then there is at least one bifurcating $L_c$ in the region surrounded by $L^{2(m-\ell-1)}_c$ and $L^{2(m-\ell)}_s$ via \eqref{outter-boun} and \eqref{inner-boun2}. It is not hard to check that there are at least $(m-2\ell-1)/2+\ell-1+1=(m-1)/2$ bifurcating $L_c$
 besides the $(m+3)/2-\ell$ bifurcating $L_c$ in \eqref{Loca-3}. Therefore, there are at least $(m-1)/2+(m+3)/2-\ell=m+1-\ell$ bifurcating $L_c$ and exactly $\ell$ bifurcating $L_s$.

In three cases above, it is not hard to check that there exist some constants $\beta_n>0$ such that
\begin{equation*}
|k_i-k_{i+1}|=|A_i-A_{i+1}|=\beta_n\delta,~i=1,...,2m-1
\end{equation*}
by \eqref{def-lam}, \eqref{Rela-del} and $k_{i+2d+1}=O(\delta^m)=o(\delta^4),~i=1,...,m-1$ by \eqref{Theta-Esti}, which implies that Lemma~\ref{lm1} holds. Therefore, analysis is under the framework of bifurcation.

{\it Step 3. Consider the case that $m=3$}. We take $\lambda_3=0$ and $\lambda_1<\lambda_2<0$ for the lower subsystem of \eqref{trans-sys}. Clearly, points $(\lambda_1,0)$ and $(0,0)$ are visible tangent points and $(\lambda_2,0)$ is an invisible tangent point. Further, orbit $\widehat\gamma^-(t,\lambda_2,0)$ intersects $\Sigma$ at $(P_1,0)$ for forward direction and at $(P_2,0)$ for backward direction. Then orbit $\widehat\gamma^+(t,P_2,0)$ intersects $\Sigma$ at $(P^+_2,0)$ for backward direction. Meanwhile, orbit $\widehat\gamma^-(t,p_1,0)$ satisfying that $\lambda_2\le p_1<0$ or $p_1\le P_1$ intersects $\Sigma$ at $(p_2,0)$ for backward direction. Then orbit $\widehat\gamma^+(t,p_2,0)$ intersects $\Sigma$ at $(p^+_2,0)$ for backward direction.

In order to investigate relationships between these points, we consider the following Cauchy Problem
\begin{equation*}
\frac{dy}{dx}=\frac{\phi(x,y)}{f^-(x,y)}(x-\lambda_1)(x-\lambda_2)x,~~~~~~y(p_1)=0.
\end{equation*}
Let $\Theta(x)$ be the solution and take $\lambda_1=-r_1\delta, \lambda_2=-\delta, p_1=-r_2\delta$, we obtain
\begin{equation*}
\Theta(p_2)=\int_{-r_2\delta}^{p_2(-r_2\delta)}\frac{\phi(s,\Theta(s))}{f^-(s,\Theta(s))}(s+r_1\delta)(s+\delta)sds=0.
\end{equation*}
By taking derivatives with respect of $\delta$ for both sides, we get
\begin{equation}
p'_2=\frac{-r^2_2(r_2-1)(r_2-r_1)R(-r_2\delta)\delta^3+\int_{-r_2\delta}^{p_2}R(s)(s(s+r_1\delta)+r_1s(s+\delta))ds}{r_2R(p_2)p_2(p_2+\delta)(p_2+r_1\delta)},
\label{x1-deri}
\end{equation}
where
\begin{equation*}
R(x):=\frac{\phi(x,\Theta(x))}{f^-(x,\Theta(x))}.
\end{equation*}
Then by taking $\delta\to0$ for both sides of \eqref{x1-deri}, we obtain
\begin{equation}
\kappa^2(r_2)h(\kappa(r_2)r_2)=h(r_2),
\label{rela-pr2}
\end{equation}
where
\begin{equation*}
\kappa(r_2):=p'_2(0),~~~h(s):=3s^2-4(1+r_1)s+6r_1.
\end{equation*}
In order to investigate $P_{1,2}$, we take $r_2=1$ and write $\eqref{rela-pr2}$ as
\begin{equation}
(\kappa(1)-1)^2(3\kappa^2(1)+2(1-2r_1)\kappa(1)+1-2r_1)=0.
\label{rela-P12}
\end{equation}
It is not hard to obtain that roots of \eqref{rela-P12} are
\begin{equation*}
\kappa_1(1)=\frac{2r_1-1+\sqrt{4r^2_1+2r_1-2}}{3},~~\kappa_2(1)=\frac{2r_1-1-\sqrt{4r^2_1+2r_1-2}}{3},~~\kappa_{3,4}(1)=1,
\end{equation*}
which implies that $P_{1,2}=\kappa_{1,2}(1)(-\delta)+O(\delta^2)$.
Then it is not hard to obtain that there exists some $r^*_1$ and $\delta_1$ such that
\begin{equation*}
\frac{|P_1|}{P_2}-\frac{|P^+_2|}{P_2}=\frac{2r_1-1+\sqrt{4r^2_1+2r_1-2}+O(\delta)}{1-2r_1+\sqrt{4r^2_1+2r_1-2}+O(\delta)}-\sqrt{\frac{V^+_{22}+O(P_2)}{V^+_{12}+O(P^+_2)}}>0
\end{equation*}
for $r_1=r^*_1$ and $\delta\in(0,\delta_1)$, i.e., $P_1<P^+_2$. Then for such $r^*_1$, it is not hard to obtain that $|\kappa(r_2)|\to1$ as $r_2\to+\infty$ by analyzing \eqref{rela-pr2} as $r_2\to+\infty$. Further, there exists some $r_{21}$ and $\delta_2$ such that
\begin{equation*}
\frac{|p_{1}|}{p_{2}}-\frac{|p^+_{2}|}{p_{2}}=\frac{r_2\delta}{|\kappa(r_2)|r_2\delta+O(\delta^2)}-\sqrt{\frac{V^+_{22}+O(p_{2})}{V^+_{12}+O(p^+_{2})}}<0
\end{equation*}
for $p_1=-r_{21}\delta$ and $\delta\in(0,\delta_2)$. To avoid unnecessary misunderstanding, we let $p_{11}:=-r_{21}\delta$ and $p^+_{21}$ denote corresponding $p^+_2$ and then obtain $p_{11}>p^+_{21}$. On the other hand, since $|\kappa(r_2)|\to1$ as $r_2\to0$, it can be proved by a similar way that there exists some $r_{22}$ such that $p_{12}:=-r_{22}\delta$ and corresponding $p^+_{22}$ satisfies that $p_{12}>p^+_{22}$. Thus, there exists $r^*_1, r_{21}, r_{22}$ and $\delta_0=\min\{\delta_1,\delta_2\}$ such that for $r_1=r^*_1$ and $\delta\in(0,\delta_0)$
\begin{equation}
P_1<P^+_2,~~~~~~p_{11}>p^+_{21},~~~~~~p_{12}>p^+_{22}.
\label{PB-3}
\end{equation}
Then it is proved in \cite{Li20} that there exists sufficiently small $\alpha\ne0$ and $\lambda_3\ne0$
such that there are two bifurcating $L_c$ appearing in a neighborhood of grazing loop
connecting $(0,0)$. Since variations of $\alpha$ and $\lambda_3$ are sufficiently small,
inequalities~\eqref{PB-3} hold and consequently, there are at least either another two
bifurcating $L_c$ if $P^+_2>\lambda_2$ or one bifurcating $L_s$ and
another one bifurcating $L_c$ if $P^+_2<\lambda_2$. Then we
take ${\boldsymbol k}={\boldsymbol 0}$, i.e., $\psi(x,{\boldsymbol 0})\equiv0$,
in system~\eqref{pws2-unfold} and obtain that there are either at least
 three bifurcating $L_c$ and exactly one bifurcating $L_s$ as shown
in Figure~\ref{Fig-SI-3}(a) or at least four bifurcating $L_c$ and no bifurcating $L_s$ as shown in Figure~\ref{Fig-SI-3}(b).
The proof of Theorem~\ref{thm1}(a) is finished for the case that $L^{gra}_*$ is of type $S$-$I$.

\begin{figure}[h]
\centering
\subfigure[3 bifurcating $L_c$, 1 bifurcating $L_s$]
 {
  \scalebox{0.38}[0.38]{
   \includegraphics{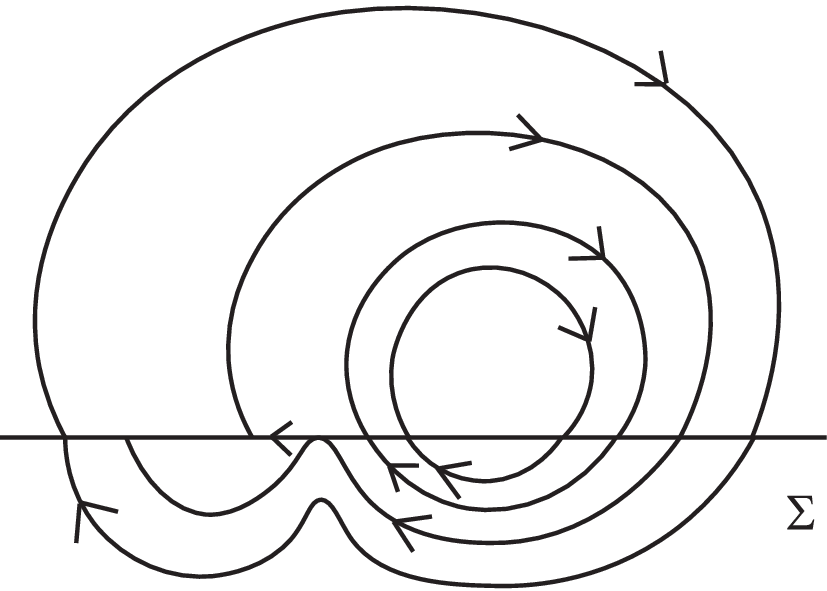}}}
\subfigure[4 bifurcating $L_c$]
 {
  \scalebox{0.38}[0.38]{
   \includegraphics{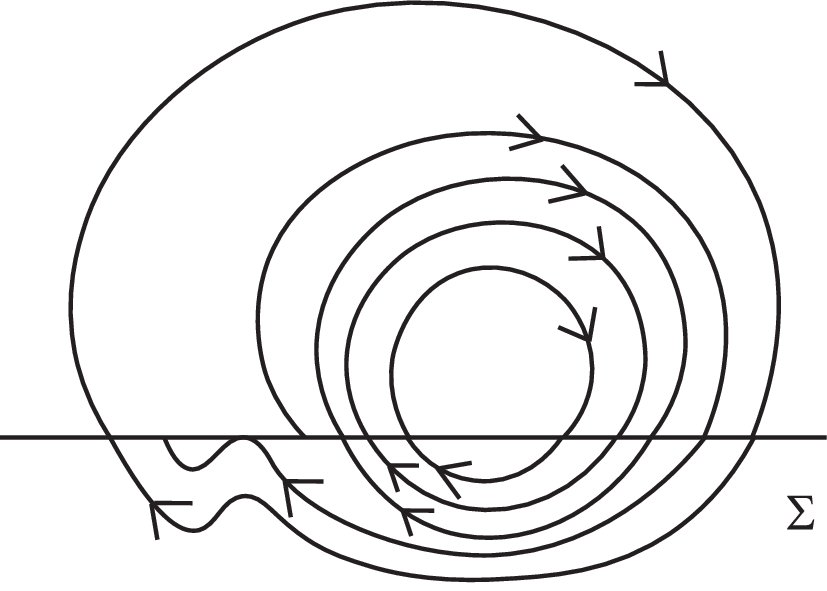}}}
   \caption{bifurcating $L_c$ and $L_s$ for $m=3$}
\label{Fig-SI-3}
\end{figure}

Now we consider the other two types for $L^{gra}_*$ stated in Theorem~\ref{thm1}(a). For loop $L^{gra}_*$ of type $S$-$L$, we split analysis into two cases: $m\ge4$ and $m=2$, which correspond to similar analysis of $m\ge 5$ and $m=3$ for type $S$-$I$ respectively. For $L^{gra}_*$ of type $S$-$V$, it is not hard to find a suitable perturbation for the lower subsystem such that
tangent point $O$ of multiplicity $(1,m)$ breaks into a left tangent point $O$ of multiplicity $(1,m-1)$ and a visible
tangent point of multiplicity $(0,1)$. Thus, this $L^{gra}_*$ becomes a grazing loop of type $S$-$L$ connecting tangent point $O$ of multiplicity $(1,m-1)$ and further bifurcation results can be
similarly obtained by analyzing the cases $m\ge5$ and $m=3$. Theorem~\ref{thm1}(a) is proved.
\end{proof}

\begin{proof}[Proof of Theorem~\ref{thm1}(b)]
For $L^{gra}_*$ of type $S$-$R$, we split the analysis into two cases: $m=2$ and $m\ge4$ like before.
In the case that $m=2$, it is not hard to find a suitable perturbation for the lower subsystem such that tangent point $O$
of multiplicity $(1,2)$ breaks into an invisible tangent point $O$ of multiplicity $(1,1)$ and
a visible tangent point of multiplicity $(0,1)$. Thus loop $L^{gra}_*$ becomes a grazing loop of type $S$-$I$ connecting one invisible tangent point $O$ of multiplicity $(1,1)$. Further, bifurcation results can be obtained by \cite{Li20}, i.e., there are at least either two bifurcating $L_c$ or one bifurcating $L_c$ and one bifurcating $L_s$.

In the case that $m\ge4$, similarly to the analysis of type $S$-$I$ we can take suitable ${\boldsymbol \lambda}$ and $\psi(x,{\boldsymbol k})$ in system~\eqref{pws2-tran} such that there exists a grazing loop connecting one $VI$ tangent point $O$ of multiplicity $(1,1)$ and $m/2$ critical loops, each of which connects a visible tangent point of multiplicity $(0,1)$.
Location of all these $1+m/2$ loops is like \eqref{Loca} and all $m/2$ critical loops
are stable. Further,  by the same method as in
the analysis in the proof of Theorem~\ref{thm1}(a) for type $S$-$I$, for each $\ell\in\{0,1,...,m/2\}$ there exists $\alpha, {\boldsymbol \lambda}, \psi(x,{\boldsymbol k})$ such that system~\eqref{pws2-unfold} has at least $m+2-\ell$ bifurcating $L_c$ and exactly $\ell$ bifurcating $L_s$.
Theorem~\ref{thm1}(b) is proved.
\end{proof}

%%%%%%%%%%%%%%%%%%%%%%%%%%%%%%%%%%%%%%%%%%%%%%%%%%%%%%%%%%%
\section{Conclusions and method discussions}
\setcounter{equation}{0}
\setcounter{lm}{0}
\setcounter{thm}{0}
\setcounter{rmk}{0}
\setcounter{df}{0}
\setcounter{cor}{0}

In this paper, we establish a relationship between $\beta_c, \beta_s$ (the numbers of bifurcating
crossing limit cycle $L_c$ and sliding loop $L_s$)
and multiplicity $m$, visibility of the tangent point on the grazing loop $L^{gra}_*$, i.e., Theorem~\ref{thm1}.
For $L^{gra}_*$ connecting one $VI$ or $VV$ tangent point,
Theorem~\ref{thm1} generalizes previous results given in \cite{Li20,Han13} for multiplicity $m$
from $m=1$ to general $m$, i.e., inequalities~\eqref{res-1} and \eqref{res-2} hold for $m\ge1$.
For $L^{gra}_*$ connecting one $VL$ or $VR$ tangent point, there is no corresponding results in previous publications and thus,
Theorem~\ref{thm1} is new.

By the idea of Poincar\'e, the method to investigate $\beta_c$ and $\beta_s$ is defining a return map
and then analyzing its fixed points. Naturally,  finding more bifurcating $L_c$ and $L_s$
is changed into finding more fixed points as in \cite{Li20,Han13} for $m=1$.
Thus, it is very important to find an appropriate domain of the return map for perturbation systems
and to find fixed points as more as possible.
However, in the case that $2\le m\le 3$, a
transition map can be defined under perturbations as in section 2 but can not be analyzed via Lemma~\ref{lm4},
which leads to a great difficulty in analyzing fixed points of the return map.
The main reason is that Lemma~\ref{lm4} relies on a segment of a known orbit but dynamical behaviors of near
orbits are unknown under perturbations. In the case that $m\ge4$, it is proved in \cite{Fang23} that the bifurcations
of a tangent point become extremely complicated as its multiplicity $m$ increasing. In details, tangent point $O$ may break into several bifurcating tangent points with low multiplicities, implying that an orbit may be tangent with $\Sigma$ for several times.
In other words, for any near orbit it is difficult to figure out the intersections between such orbit and $\Sigma$ under perturbations,
i.e., the domain of the return map is unclear.
Therefore, there are essential difficulties in the analysis of $m\ge 2$ compared with $m=1$.

To overcome the above difficulties, we construct parameter $\alpha\in \mathbb{R}$ to control intersections
between the limit cycle of the upper subsystem and $\Sigma$, and construct parameter ${\boldsymbol \lambda}\in \mathbb{R}^m$
to control locations and multiplicities of bifurcating tangent points. In the case that $m\ge4$,
function $\psi(x,{\boldsymbol k})$ is defined to control intersections between $\Sigma$ and
tangent orbits. Further, dynamical behaviors with respect to return are clear under suitable
perturbations. Results \eqref{res-2} and \eqref{res-4} are obtained via constructing different
$\psi(x,{\boldsymbol k})$.
In the case that $2\le m\le3$, by Lemma~\ref{lm1} function $\psi(x,{\boldsymbol k})$ is restricted to be
zero function to ensure that \eqref{pws2-unfold} is a small perturbation of \eqref{pws2}.
This restriction leads to that we have to find the information of the transition map for the lower subsystem
but, the information can not be obtained by Lemma~\ref{lm4} because of the lack of one known orbit of the lower subsystem.
In our method, we unfold the lower
subsystem along some specific curves in the parameter space, which helps give the information of the transition map
and then the dynamical behaviors with respect to return.

We have to remark that system~\eqref{pws2} is a specific piecewise-smooth system having a tangent point of multiplicity $(1,m)$
because $g^-(x,y)$ is restricted to be $\phi(x,y)x^m$.
This restriction helps us analyze the perturbation system~\eqref{pws2-unfold} and construct $\psi(x,{\boldsymbol k})$
to obtain bifurcating $L_c$ and $L_s$.
However, for general $g^-(x,y)$ our method in this paper is invalid and it is under our futural consideration.

%%%%%%%%%%%%%%%%%%%%%%%%%%%%%%%%%%%%%%%%%%%%%%%%%%%%%%%%%%%%%%%%%


{\footnotesize

\begin{thebibliography}{99}



\bibitem{Bonet18}
C. Bonet-Reves, J. Larrosa, T. M. Seara,
Regularization around a generic codimension one fold-fold singularity,
{\it J. Differential Equa.} {\bf 265}(2018), 1761-1838.

\bibitem{Bernardo08}
M. di Bernardo, C. J. Budd, A. R. Champneys, P. Kowalczyk,
{\it Piecewise-Smooth Dynamical Systems: Theory and Applications},
Springer-Verlag, London, 2008.

\bibitem{Bernardo082}
M. di Bernardo, C. J. Budd, A. R. Champneys, P. Kowalczyk, A. Nordmark, G. Tost, P. Piiroinen,
Bifurcations on nonsmooth dynamical systems,
{\it SIAM Rev.} {\bf 50}(2008), 629-701.

\bibitem{Ponce22}
M. Esteban, E. Freire, E. Ponce, F. Torres,
On normal forms and return maps for pseudo-focus points,
{\it J. Math. Anal. Appl.} {\bf 507}(2022), Paper No. 125774.

\bibitem{Fang21}
Z. Fang, X. Chen,
Global dynamics of a piecewise smooth system with a fold-cusp and general parameters,
{\it Qual. Theo. Dyn. Syst.} {\bf 21}(2022), Paper No. 55.

\bibitem{Fang23}
Z. Fang, X. Chen,
Classifications and bifurcations of tangent points and their loops of planar piecewise-smooth systems,
in preprint.

\bibitem{Filippov88}
A. F. Filippov,
{\it Differential Equations with Discontinuous Righthand Sides},
Kluwer Academic Publishers, Dordrecht, 1988.

\bibitem{Ponce14}
E. Freire, E. Ponce, F. Torres,
A general mechanism to generate three limit cycles in planar Filippov systems with two zones,
{\it Nonlinear Dyn.} {\bf 78}(2014), 251-263.

\bibitem{Ponce15}
E. Freire, E. Ponce, F. Torres,
On the critical crossing cycle bifurcation in planar Filippov systems,
{\it J. Differential Equa.} {\bf 259}(2015), 7086-7107.

\bibitem{Teixeira11}
M. Guardia, T. M. Seara, M. A. Teixeira,
Generic bifurcations of low codimension of planar Filippov systems,
{\it J. Differential Equa.} {\bf 250}(2011), 1967-2023.

\bibitem{JKHaleBook}
J. K. Hale,
{\it Dynamics and Bifurcations},
Springer-Verlag, New York, 1991.

\bibitem{Kuznetsov03}
Yu. A. Kuznetsov, S. Rinaldi, A. Gragnani,
One-parameter bifurcations in planar Filippov systems,
{\it Int. J. Bifurc. Chaos} {\bf 13}(2003), 2157-2188.

\bibitem{KuznetsovBook}
Yu. A. Kuznetsov,
{\it Elements of Applied Bifurcation Theory},
Springer-Verlag, New York, 2004.

\bibitem{Lee18}
J. M. Lee,
{\it Introduction of Smooth Manifold},
Springer-Verlag, New York, 2012.

\bibitem{Li20}
T. Li, X. Chen,
Degenerate grazing-sliding bifurcations in planar Filippov systems,
{\it J. Differential Equa.} {\bf 269}(2020), 11396-11434.

\bibitem{Han12}
F. Liang, M. Han,
Degenerate Hopf bifurcation in nonsmooth planar systems,
{\it Int. J. Bifurc. Chaos} {\bf 22}(2012), Paper No. 1250057.

\bibitem{Han13}
F. Liang, M. Han, X. Zhang,
Bifurcation of limit cycles from generalized homoclinic loops in planar piecewise smooth systems,
{\it J. Differential Equa.} {\bf 255}(2013), 4403-4436.

\bibitem{Novaes18}
D. D. Novaes, M. A. Teixeira, I. O. Zeli,
The generic unfolding of a codimension-two connection to a two-fold singularity of planar Filippov systems,
{\it Nonlinearity} {\bf 31}(2018), 2083-2104.

\bibitem{Siller}
T. E. Siller,
A missing generic local fold-fold bifurcation in planar Filippov systems,
{\it Int. J. Bifurc. Chaos} {\bf 32}(2022), Paper No. 2250031.

\bibitem{Huang22}
F. Wu, L. Huang, J. Wang,
Bifurcation of the critical crossing cycle in a planar piecewise smooth system with two zones,
{\it Disc. Cont. Dyn. Syst. Series B} {\bf 27}(2022), 5047-5083.

\bibitem{ZhangJ98}
J. Zhang, B. Feng,
{\it Geometric Theory and Bifurcation Problems in Ordinary Differential Equations},
Peking University Press, Beijing, 1981 (in Chinese).




\end{thebibliography}
}
\end{document}